\documentclass[a4paper,12pt]{amsart}
\usepackage{amsmath,amssymb}
\usepackage{graphicx}
\newtheorem{thm}{Theorem}
\newtheorem{prop}[thm]{Proposition}
\newtheorem{cor}[thm]{Corollary}

\begin{document}

\title{Areas associated to a quadrilateral}

\author{Oleg Mushkarov}
\address{O. Mushkarov\\Institute of Mathematics and Informatics\\Bulgarian Academy
of Sciences\\Acad. G. Bonchev 8, 1113 Sofia, Bulgaria}

\email{muskarov@math.bas.bg}

\author{Nikolai Nikolov}
\address{N. Nikolov\\Institute of Mathematics and Informatics\\Bulgarian Academy
of Sciences\\Acad. G. Bonchev 8, 1113 Sofia, Bulgaria
\vspace{1mm}
\newline Faculty of Information Sciences\\
State University of Library Studies and Information Technologies\\
Shipchenski prohod 69A, 1574 Sofia,
Bulgaria}

\email{nik@math.bas.bg}


\maketitle

\begin{abstract} We study the relationship between the areas of the consecutive
quadrilaterals cut from a convex quadrilateral in the plane by
means of a finite or infinite number of straight lines
intersecting two of its opposite sides. Moreover, we obtain a
geometric description of all possible areas obtained in this way
given the ratios of the lengths of consecutive segments the lines
divide these two opposite sides.
\end{abstract}

\section{Introduction}
Let $(p_1,\dots,p_n)$ and $(p_1',\dots,p_n')$ be given $n$-tuples
of positive real numbers. For a convex quadrilateral $ABCD$ in the
plane denote by $A_1,A_2,\dots,A_{n-1}$ and $D_1,D_2,\dots,$
$D_{n-1}$ the points on the sides $AB$ and $CD$, respectively such
that $AA_1:A_1A_2:\dots:A_{n-1}B=p_1:p_2:\dots:p_n$ and
$DD_1:D_1D_2:\dots : D_{n-1}C=p_1':p_2':\dots:p_n'$ (Fig. 1). The
purpose of this note is to find a geometric description of the set
$\mathcal{A}$ of all possible $n$-tuples of areas
$(\mathcal{A}_1,\mathcal{A}_2,\dots,\mathcal{A}_n)$, where
$\mathcal{A}_1= [AA_1D_1D], \mathcal{A}_2= [A_1A_2D_2D_1], \dots ,
\mathcal{A}_n=[A_{n-1}BCD_{n-1}]$, when $ABCD$ runs over all
convex quadrilaterals in the plane. This problem is motivated by
the linear

\begin{center}
\setlength{\unitlength}{20pt} \linethickness{0.5pt}
\begin{picture}(12,7)
\put(0,0){\line(4,3){0.5}}
    \put(0.8,0.6){\line(4,3){0.5}}
    \put(1.6,1.2){\line(4,3){0.5}}
    \put(2.4,1.8){\line(4,3){0.5}}
    \put(3.2,2.4){\line(4,3){0.5}}
    \put(4,3){\line(4,3){0.5}}
\put(0,0){\line(1,0){0.6}}
    \put(1,0){\line(1,0){0.6}}
    \put(2,0){\line(1,0){0.6}}
    \put(3,0){\line(1,0){0.6}}
\put(3.8,0){\line(1,4){0.85}} \put(4.65,3.45){\line(4,3){4.1}}
\put(3.8,0){\line(1,0){8.2}} \put(7,0){\line(-1,4){1.1}}
\put(9,0){\line(-1,3){1.8}} \put(12,0){\line(-1,2){3.25}}

\put(0,0){\circle*{0.2}} \put(-0.2,-0.6){$E$}
\put(3.8,0){\circle*{0.2}} \put(3.6,-0.6){$A$}
\put(4.65,3.45){\circle*{0.2}} \put(4.4,3.7){$D$}
\put(7,0){\circle*{0.2}} \put(6.8,-0.6){$A_1$}
\put(9,0){\circle*{0.2}} \put(8.8,-0.6){$A_{n-1}$}
\put(12,0){\circle*{0.2}} \put(11.8,-0.6){$B$}
\put(8.75,6.5){\circle*{0.2}} \put(8.5,6.7){$C$}
\put(7.2,5.35){\circle*{0.2}} \put(6.35,5.75){$D_{n-1}$}
\put(5.9,4.4){\circle*{0.2}} \put(5.5,4.7){$D_1$}
\put(3,1.3){$\mathcal{A}_0$} \put(5.2,1.7){$\mathcal{A}_1$}
 \put(9.2,2.3){$\mathcal{A}_n$}
\end{picture}

\vspace*{1cm}

Fig. 1

\end{center}

\noindent and analytic relations between $\mathcal{A}_1,$
$\mathcal{A}_2$ and $\mathcal{A}_3$ in the case when $n=3$ and
$p_1:p_2:p_3=p_1':p_2':p_3'$ found in \cite{NN}. Here we consider
the case of arbitrary $n$ and arbitrary ratios (Section 2 and
Section 3) as well as that when the sides of the quadrilaterals
are divided into infinitely many segments (Section 4). We also
show that the last case can be reduced to that when $n=3$ and
$n=2$ (Section 5).

Understanding the article requires a basic knowledge of linear
algebra and analytic geometry

\section{Area attainable points in $\Bbb R^3$} We will first consider the case when
the sides of the quadrilaterals are divided into $3$ segments. Let
their ratios be determined by the triples $(p_1,p_2,p_3)$ and
$(p_1',p_2',p_3')$.
\smallskip

{\bf Definition.} A point $(x_1,x_2,x_3)$ in the $3$-dimensional
Euclidean space $\Bbb R^3$ is called {\it area attainable} if
there is a convex quadrilateral $ABCD$ in the plane such that
$(\mathcal{A}_1,\mathcal{A}_2,\mathcal{A}_3)=(x_1,x_2,x_3).$
\smallskip

It is clear that the set $\mathcal{A}$ of all area attainable
points is the union of rays in the open first octant $\Bbb R_+^3$
with a common initial point $O=(0,0,0)$.

\begin{prop}\label{PON1} Set
$\Delta=(p_1+p_2+p_3)p_1'p_3'p_2-(p_1'+p_2'+p_3')p_1p_3p_2'.$

a) If $ \Delta \neq 0,$ then $\mathcal{A}$ is the union of two
nonintersecting open trihedral angles with a common vertex $O$ and
a common face, and an open ray with origin $O,$ lying on this face
and different from its edges.

b) If $\Delta =0,$ then $\mathcal{A}$ is an open angle with vertex
$O$ in the plane
\begin{equation}\label{1}
\alpha:\
(p_2p_3'-p_3p_2')x_1+(p_3p_1'-p_1p_3')x_2+(p_1p_2'-p_2p_1')x_3=0
\end{equation}
when $p_1:p_1'\neq p_3:p_3'$ or in the plane
\begin{equation}\label{2}
\alpha:\
\frac{p_2+p_3}{p_1}x_1-\frac{p_1+2p_2+p_3}{p_2}x_2+\frac{p_1+p_2}{p_3}x_3=0
\end{equation}
 when $p_1:p_1'=p_3:p_3'$.
\smallskip

Note that if $\Delta=0$ and two of the ratios $p_1:p_1',$
$p_2:p_2',$ $p_3:p_3'$ are equal, then the third one is also equal
to them.
\end{prop}

\noindent {\bf Proof.} Set $\vec v_1=(p_1,p_2,p_3),$ $\vec
v_2=(p_1',p_2',p_3'),$ $\vec v_0=\vec v_1+\vec v_2,$
$$\vec v_3=(p_1p_1',p_2p_2'+p_2p_1'+p_2'p_1,p_3p_3'+p_3(p_1'+p_2')+p_3'(p_1+p_2)),$$
$$\vec v_4=(p_1p_1'+p_1(p_2'+p_3')+p_1'(p_2+p_3),p_2p_2'+p_2p_3'+p_2'p_3,p_3p_3'),$$
and denote by $l_j=\{t\vec v_j:t>0\}$ the rays determined by these
vectors.

Note first that if $AB\parallel CD$, then $l_0\subset
\mathcal{A}.$ Let now $AB\nparallel CD$ and denote by
$\mathcal{Q}_1$ the set of convex quadrilaterals $ABCD$ in the
plane such that $A$ lies between $B$ and $E=AB\cap CD$ (the set
$\mathcal{Q}_2$ is defined analogously when $B$ lies between $A$
and $E$). Set $EA:AA_1=p_0:p_1,$ $ED:DD_1=p_0':p_1',$
$\mathcal{A}_0=[EAB]=sp_0p_0'.$ Then
$$\frac{\mathcal{A}_0+\mathcal{A}_1}{\mathcal{A}_0}=\frac{(p_0+p_1)(p_0'+p_1')}{p_0p_0'},$$
$$\frac{\mathcal{A}_0+\mathcal{A}_1+\mathcal{A}_2}{\mathcal{A}_0}=\frac{(p_0+p_1+p_2)(p_0'+p_1'+p_2')}{p_0p_0'},$$
$$\frac{\mathcal{A}_0+\mathcal{A}_1+\mathcal{A}_2+\mathcal{A}_3}{S_0}=\frac{(p_0+p_1+p_2+p_3)(p_0'+p_1'+p_2'+p_3')}{p_0p_0'},$$
from where we get
$$ \mathcal{A}_1=s(p_1p_0'+p_1'p_0+p_1p_1'),$$
$$ \mathcal{A}_2=s(p_2p_0'+p_2'p_0+p_2p_2'+p_2p_1'+p_2'p_1),$$
$$ \mathcal{A}_3=s(p_3p_0'+p_3'p_0+p_3p_3'+p_3(p_1'+p_2')+p_3'(p_1+p_2)).$$

Let $M_1$ be the $3\times 3$ matrix with rows $\vec v_1,$ $\vec
v_2,$ $\vec v_3.$ Then the identities above may be written in the
following matrix form
$(\mathcal{A}_1,\mathcal{A}_2,\mathcal{A}_3)=(sp_0',sp_0,s)M_1.$
Since the map $ABCD\to (sp_0',sp_0,s)$ is a surjection from $Q_1$
onto the open first octant $\Bbb R_+^3$ of $\Bbb R^3$ (take, for
example, $\angle AED=90^\circ,$ $AE=2p_0$ and $DE=sp_0'$), it
follows that the set $\mathcal{A}'$ of area attainable points for
which $A$ lies between $E$ and $B$ is the image of $\Bbb R_+^3$
under the linear transformation $M_1$ of $\Bbb R^3$ with matrix
$M_1.$

Analogously, the set $\mathcal{A}''$ of area attainable points for
which $B$ lies between $A$ and $E$ is the image of $\Bbb R_+^3$
under the linear transformation $M_2$ of $\Bbb R^3$ with matrix
$M_2$ whose rows are  $\vec v_1,$ $\vec v_2,$ $\vec v_4.$

It is clear that the set $\mathcal{A}$ of all area attainable
points is given by $\mathcal{A}=l_0\cup \mathcal{A}'\cup
\mathcal{A}''$ and now we will describe geometrically this set
depending on whether $\Delta \neq 0$ or $\Delta=0,$ i.e. if
$\mbox{rank }M_1=\mbox{rank }M_2$ is equal to $3$ or $2$ (the
ranks of $M_1$ and $M_2$ are never equal to $1$.)

Let $\Delta \neq 0.$ Then $\mathcal{A}'$ and $\mathcal{A}''$ are
the open trihedral angles with edges the rays $l_1,l_2,l_3$ and
$l_1,l_2,l_4$ (the images of the positive axes under the linear
transformations $M_1$ and $M_2,$ respectively), and $l_0$ lies on
the face with edges $l_1$ and $l_2$ (Fig. 2). It is easy to check
that the scalar triple product of the vectors $\vec v_1,$ $\vec
v_2,$ $\vec v_3$ is equal to $\det M_1=-\Delta,$ and that of $\vec
v_1,$ $\vec v_2,$ $\vec v_4$ is equal to $\det M_2=\Delta .$ In
particular, these two triples of vectors determine opposite
orientations on $\Bbb R^3$ and therefore $\mathcal{A}'\cap
\mathcal{A}''=\varnothing.$ Note that this follows also from the
identity
$$\vec v_3+\vec v_4=(p'_1+ p'_2+p'_3)\vec v_1+ (p_1+ p_2+p_3)\vec v_2.$$

\begin{center}
\setlength{\unitlength}{20pt} 
\linethickness{0.5pt}
\begin{picture}(6,6)
\put(0,0){\vector(1,4){1.2}} \put(0,0){\line(4,3){0.5}}
    \put(0.8,0.6){\line(4,3){0.5}}
    \put(1.6,1.2){\line(4,3){0.5}}
    \put(2.4,1.8){\line(4,3){0.5}}
    \put(3.2,2.4){\line(4,3){0.5}}
    \put(4,3){\vector(4,3){0.5}}
\put(0,0){\line(4,1){0.5}}
    \put(0.8,0.2){\line(4,1){0.5}}
    \put(1.6,0.4){\line(4,1){0.5}}
    \put(2.4,0.6){\line(4,1){0.5}}
    \put(3.2,0.8){\line(4,1){0.5}}
    \put(4,1){\vector(4,1){0.5}}
\put(0,0){\vector(1,0){6}} \put(0,0){\vector(4,-1){5}}
\put(0,0){\circle*{0.2}} \put(-0.4,-0.4){$0$} \put(0.2,3){$l_4$}
\put(2,2){$l_2$} \put(3,1){$l_0$} \put(4,0.2){$l_1$}
\put(3.8,-1.5){$l_3$}
\end{picture}

\vspace*{1.2cm}

Fig. 2

\end{center}

Let now $\Delta=0.$ Then $\det M_1=\det M_2=0$ and the rays $l_0,$
$l_1,$ $l_2,$ $l_3,$ $l_4$ lie on a plane $\alpha\supset
\mathcal{A}.$ It is defined by \eqref{1} or \eqref{2}, which
follows respectively from $\alpha\perp\vec v_1\times\vec v_2$ or
$\alpha\perp\vec v_3\times\vec v_4.$

Note that the image $M_1(\Bbb R_+^3)$ (resp. $M_2(\Bbb R_+^3)$)
consists of the points in $\Bbb R^3$ whose vectors are linear
combinations of $\vec v_1, \vec v_2, \vec v_3$ (resp. $\vec v_1,
\vec v_2, \vec v_4$)  with positive coefficients. A direct check
shows that $\vec v_1$ and $\vec v_2$ are linear combinations of
$\vec v_3$ and $\vec v_4$ with positive coefficients. Hence in the
case $p_1:p_1'\neq p_3:p_3',$ we have
$\mathcal{A}'=\angle(l_1,l_3),$ $\mathcal{A}''=\angle(l_2,l_4)$ if
$l_2\subset\angle(l_1,l_3)$ and $\mathcal{A}'=\angle(l_2,l_3),$
$\mathcal{A}''=\angle(l_1,l_4)$ if $l_1\subset\angle(l_2,l_3);$ in
particular, $l_0\subset \angle(l_1,l_2)=\mathcal{A}'\cap
\mathcal{A}''.$ In the case $p_1:p_1'=p_3:p_3'$ it follows that
$\mathcal{A}'=\angle(l_0,l_3)$ and
$\mathcal{A}''=\angle(l_0,l_4).$ Hence in both cases $\mathcal
{A}=\angle(l_3,l_4).$\qed
\smallskip

The above arguments imply the following

\begin{cor}\label{CON1} If $AB\parallel CD,$ then
\begin{equation}\label{3}
\mathcal{A}_1:\mathcal{A}_2:\mathcal{A}_3=(p_1+p_1'):(p_2+p_2'):(p_3+p_3').
\end{equation}

Conversely, if $\Delta\neq 0$ and $\eqref{3}$ is fulfilled, then
$AB\parallel CD.$
\end{cor}

\noindent {\bf Remark.} a) If $p_1:p_2:p_3=p_1':p_2':p_3'$, then
$AB\parallel CD$ follows from the weaker assumptions
$\mathcal{A}_1:\mathcal{A}_2=p_1:p_2,$
$\mathcal{A}_2:\mathcal{A}_3=p_2:p_3$ or
$\mathcal{A}_3:\mathcal{A}_1=p_3:p_1.$

b) If  $\Delta=0,$ but $p_1:p_1'\neq p_3:p_3',$ then there are
quadrilaterals from $\mathcal{Q}_1$ and $\mathcal{Q}_2,$ such that
\eqref{3} is fulfilled.
\smallskip

The proof of Proposition~\ref{PON1} implies also an analytic
description of the area attainable points $\mathcal{A}.$ For
example, we have the following

\begin{cor}\label{CON2} (\cite[Proposition 3]{NN}) If
$p_1:p_2:p_3=p_1':p_2':p_3',$ then the area attainable points
$(x_1,x_2,x_3)$ are those for which \eqref{2} is fulfilled,
$x_1>0$ and
\begin{equation}\label{4}
\frac{p_3^2}{p_1(p_1+2p_2+2p_3)}<\frac{x_3}{x_1}<\frac{p_3(2p_1+2p_2+p_3)}{p_1^2}.
\end{equation}

\end{cor}

\noindent {\bf Proof.} Note that
$$l_3=\{t(p_1^2,p_2(p_2+2p_1),p_3(2p_1+2p_2+p_3)):t>0\},$$
$$l_4=\{t(p_1(p_1+2p_2+2p_3),p_2(p_2+2p_3),p_3^2):t>0\}$$
and it follows from the proof of Proposition~\ref{PON1} that the
area attainable points $(x_1,x_2,x_3)$ are those for which
\eqref{2} is fulfilled, and
$$x_1=\lambda p_1^2+\mu p_1(p_1+2p_2+2p_3),\ x_3=\lambda p_3(2p_1+2p_2+p_3)+\mu p_3^2,\ \lambda,\mu>0.$$
Solving this system with respect to $\lambda$ and $\mu$ we see
that it is equivalent to $x_1>0$ and \eqref{4}.\qed

\section{Area attainable points in $\Bbb R^n$} In this section we will obtain
a geometric description of the set $\mathcal{A}$ of area
attainable points when the sides $AB$ and $CD$ of a convex
quadrilateral $ABCD$ in the plane are divided into $n$ segments.
In this case $\mathcal{A}$ is a subset of the $n$-dimensional
Euclidean space $\Bbb R^n$.

Let $(p_1,\dots,p_n)$ and $(p_1',\dots,p_n')$ be two $n$-tuples of
positive real numbers which determine the ratios of the
consecutive $n$ segments of which the sides $AB$ and $CD$ are
divided. The case when these two $n$-tuples coincide has been
considered in \cite {NN} , where an analytic description of the
set $\mathcal{A}$ is given.

If $n=2$, then the proof of Proposition~\ref{PON1} for
$p_3=p_3'=0$ implies that $\mathcal{A}$ is an open angle with
vertex $O$ and arms
$$\{t(p_1p_1',p_2p_2'+p_2p_1'+p_2'p_1): t>0\},\quad
\{t(p_1p_1'+p_1p_2'+p_1'p_2,p_2p_2'): t>0\}.$$

The analog of Proposition~\ref{PON1} in higher dimensions is the
following

\begin{prop}\label{PON2} Let $n\ge 4$ and set
$$\Delta_i=(p_{i-1}+p_i+p_{i+1})p_{i-1}'p_{i+1}'p_i-(p_{i-1}'+p_i'+p_{i+1}')p_{i-1}p_{i+1}p_i',
\ 2\leq i \leq n-1.$$ a) If at least one of the numbers $\Delta_i$
is different from $0,$ then $\mathcal{A}$ is the union of two
nonintersecting open trihedral angles in the same $3$-dimensional
vector subspace of  $\Bbb R^n$ with a common vertex
$O=(0,0,\dots,0)$ and a common face, and an open ray with origin
$O$ lying in this face and different from its edges;

\noindent b) If all numbers $\Delta_i$ are equal to $0,$ then
$\mathcal{A}$ is an open angle with vertex $O.$

\end{prop}

\noindent {\bf Sketch of proof.} Set $\vec v_1=(p_1,\dots,p_n),$
$\vec v_2=(p_1',\dots,p_n'),$ $\vec v_0=\vec v_1+\vec v_2,$ $\vec
v_3=(s_1,\dots,s_n),$ $\vec v_4=(t_1,\dots,t_n),$ where
$$s_i=-p_ip_i'+p_i\sum_{j=1}^ip_j'+p_i'\sum_{j=1}^ip_j,\quad
t_i=-p_ip_i'+p_i\sum_{j=i}^np_j'+p_i'\sum_{j=i}^np_j,\quad 1\leq
i\leq n.$$ Using the above notations and those in Section 1 we see
that if $AB\nparallel CD$, then the areas
$\mathcal{A}_1,\mathcal{A}_2,\dots,\mathcal{A}_n$ are given by the
following formulas:
$$ \mathcal{A}_i=s(p_0p'_i+p_0'p_i+s_i),\quad 1\leq i\leq n $$
or
$$ \mathcal{A}_i=s(p_0p'_i+p_0'p_i+t_i),\quad 1\leq i\leq n $$
depending on whether $A$ lies between $E=AB\cap CD$ and $B$ or $B$
lies between $A$ and $E$.

Denote by $l_j=\{t\vec v_j:t>0\}$ the rays determined by the
vectors $v_j, 0\leq j \leq 4$.

In case b) the set $\mathcal{A}$ is an open angle with vertex $O$
and arms $l_3$ and $l_4,$ which determine a plane $\alpha$ in
$\Bbb R^n.$ Let us note that in this case either $p_j:p_k\neq
p_j':p_k'$ for $j\neq k,$ or
$p_1:p_2:\dots:p_n=p_1':p_2':\dots:p_n'.$ To see this we set
$q_i=p_i/p_{i-1},$ $q_i'=p_i'/p_{i-1}'$ and rewrite the identities
$\Delta_i=0$ in the form
$$q_{i+1}-q'_{i+1}=\frac{(1+q_{i+1}')q_{i+1}}{(1+q_i)}(q_i-q_i').$$
It follows that if $q_2\neq q'_2$, then the sequence with general
term $p_{i+1}/p'_{i+1}$ us strictly monotonic and therefore
$p_j:p_k\neq p_j':p_k'$ for $j\neq k$. If $q_2= q'_2,$ then
$q_{i+1}=q'_{i+1}, $ i.e.
$p_1:p_2:\dots:p_n=p_1':p_2':\dots:p_n'.$

Let us note also that the plane $\alpha$ is the intersection of
the linearly independent hyperplanes $\alpha_i, 2\leq i \leq n-1$,
defined by the equations:
\begin{equation}\label{5}
\alpha_i:(p_ip_{i+1}'-p_{i+1}p_i')x_{i-1}+(p_{i+1}p_{i-1}'-p_{i-1}p_{i+1}')x_i+
(p_{i-1}p_i'-p_ip_{i-1}')x_{i+1}=0
\end{equation}
in case a) and
\begin{equation}\label{6}\alpha_i:\frac{p_{i+1}+p_i}{p_{i-1}}x_{i-1}-
\frac{p_{i-1}+2p_i+p_{i+1}}{p_i}x_i+\frac{p_{i-1}+p_i}{p_{i+1}}x_{i+1}=0
\end{equation}
in case b).

In case a) we have $\mathcal{A}=l_0\cup \mathcal{A}'\cup
\mathcal{A}'',$ where $\mathcal{A}'$ is a trihedral angle with
edges $l_1,l_2,l_3$, and $\mathcal{A}''$ -- that with edges
$l_1,l_2,l_4.$ The identity
$$\vec v_3+\vec v_4=\vec v_1\sum_{j=1}^np_j'+\vec v_2\sum_{j=1}^np_j$$
shows that these two trihedral angles lie in the same
$3$-dimensional vector subspace of $\Bbb R^n$ and $l_3$, and $l_4$
lie on different sides of the plane determined by $l_1$ and $l_2.$
This $3$-dimensional vector subspace can be described as the
intersection of $n-3$ linearly independent hyperplanes $\beta_i.$
For example, if $\Delta_k\neq 0 $ for some $k$ we may choose
\begin{equation}\label{9}\beta_i:\ x_i+c_{k-1,i}x_{k-1}+c_{k,i}x_k+
c_{k+1,i}x_{k+1}=0,\ i=1,\dots,k-2,k+2,\dots,n,
\end{equation}
where the coefficients $c_{k-1,i},c_{k,i}, c_{k+1,i}$ are the
solutions $x,y,z$ of the system
$$p_{k-1}x+p_ky+p_{k+1}z=-p_i$$
$$p'_{k-1}x+ p'_ky+p'_{k+1}z=-p'_i$$
$$s_{k-1}x+s_ky+s_{k+1}z=-s_i.$$
The explicit formulas for these solutions are given by Kramer's
formulas.\qed

\section{Area attainable points in $\ell^1$}

Now we will consider the case when the sides of the quadrilaterals
are divided into infinitely many segments.

Let $p=(p_i)_{i\in\Bbb N}$ and $p'=(p_i')_{i\in\Bbb N}$ be two
infinite sequences of positive integers. We will say that an
infinite sequence $x=(x_i)_{i\in\Bbb N}$ of real numbers is area
attainable if there is a quadrilateral $ABCD$ in the plane and
points $A_0=A, A_1, A_2,\dots$ and $D_0=D, D_1, D_2,\dots $ on the
sides $AB$ and $DC$, such that $A_0A_1:A_1A_2:\dots=p_1:p_2:\dots
$, $D_0D_1:D_1D_2:\dots=p_1':p_2':\dots$ and
$x_i=\mathcal{A}_i=[A_{i-1}A_iD_iD_{i-1}]$ for all $i\in\Bbb N.$

Denote by  $\ell^1$ the space of all sequences $r=(r_i)_{i\in\Bbb
N}$ of real numbers such that $\sum_{i=1}^\infty|r_i|<\infty.$ It
is clear that the set $\mathcal{A}$ of area attainable points is
empty if $p\not\in \ell^1$ or $p'\not\in \ell^1.$

\begin{prop}\label{PON3} Let $p,p'\in \ell^1$ be given infinite sequences
of positive real numbers and let $\Delta_i$ ($i\ge 2$) be the
numbers defined in Proposition~\ref{PON2}.

a) If at least one of the numbers $\Delta_i$ is different from
$0,$ then $\mathcal{A}$ is the union of two nonintersecting open
trihedral angles in the same $3$-dimensional vector subspace of
$\ell^1$ with a common vertex $O=(0,0,\dots,0)$ and a common face,
and an open ray with origin $O$ lying on this face and different
from its edges.

b) If all numbers $\Delta_i$ are equal to $0,$ then $\mathcal{A}$
is an open ($2$-dimensional) angle in $\ell^1$ with vertex $O.$

\end{prop}

The proof of this proposition is the same as that of
Proposition~\ref{PON2} with some obvious changes. For example, in
the description of the set $\mathcal{A}''$ we may assume that
$A_n\to B,$ $D_n\to C$ and then the identity
$$\frac{\sum_{j=i}^\infty
S_j}{S_0}=\frac{\left(p_0+\sum_{j=i}^\infty p_j\right)\left(p_0'+
\sum_{j=i}^\infty p_j'\right)}{p_0p_0'}$$ implies that in this
case
$$ \mathcal{A}_i=s(p_0p'_i+p_0'p_i+t_i),$$
where
$$t_i=-p_ip_i'+p_i\sum_{j=i}^\infty p_j'+p_i'\sum_{j=i}^\infty p_j,\quad i\in\Bbb N.$$
Then $\vec v_3, \vec v_4\in \ell^1$ since
$$\sum_{i=1}^\infty(s_i+t_i)=2\sum_{i=1}^\infty p_i\sum_{i=1}^\infty p_i'.$$

The three  and two dimensional vector subspaces of $\ell^1$
containing the set $\mathcal{A}$ of area attainable sequences can
be defined as intersections of countable many linearly independent
hyperplanes in $\ell^1$. Note also that Proposition~\ref{PON2} is
a particular case of Proposition~\ref{PON3} in view of the natural
embeding of $\Bbb R^n$ in $\ell^1.$

\section{Back to $\Bbb R^3$ and $\Bbb R^2.$}

The considerations in previous sections allow us to reduce the
case of $\ell^1$ (in particular, $\Bbb R^n$) to the cases $\Bbb
R^3$ or $\Bbb R^2.$

\begin{prop}\label{PON4} Let $p,p'\in \ell^1.$

a) If $\Delta_k\neq 0$ for some $k\ge 2,$ then $x\in \mathcal{A}'$
precisely when \eqref{9} is fulfilled for all $i\neq k-1, k, k+1$ and
$(\sum_{j=1}^{k-1}x_j, x_k,x_{k+1})$ is an area attainable point in $\Bbb R^3$ 
corresponding to the triples
$(\sum_{j=1}^{k-1}p_j,p_k,p_{k+1})$ and
$(\sum_{j=1}^{k-1}p'_j,p'_k,p'_{k+1}),$ 
and such that $A$ is between $B$ and $E.$

Analogously, $x\in {\mathcal{A}}''$ precisely when \eqref{9} is
fulfilled for all $i\neq k-1, k, k+1$ and $(x_{k-1},x_k,\sum_{j=k+1}^\infty
x_j)$ is an area attainable point in $\Bbb R^3$ corresponding to the triples
$(p_{k-1},p_k,\sum_{j=k+1}^\infty p_j)$ and
$(p_{k-1}',p_k',\sum_{j=k+1}^\infty p_j'),$ 
and such that $B$ is between $A$ and $E.$

b) If $\Delta_k=0$ for all $k\ge 2,$ then $x\in\mathcal{A}$
precisely when $x_1>0,$ \eqref{5} or \eqref{6} is fulfilled for
all $i\ge 2$ and
$$\frac{\Sigma_2}{\Sigma_1}<\frac{x_2}{x_1}<\frac{p_2p_2'+p_2p_1'+p_2'p_1}{p_1p_1'},$$
where
$$\Sigma_1=p_1p_1'+p_1\sum_{j=2}^\infty p_j'+p_1'\sum_{j=2}^\infty p_j,\
\Sigma_2=p_2p_2'+p_2\sum_{j=3}^\infty p_j'+p_2'\sum_{j=3}^\infty
p_j$$

\end{prop}

We leave the proof to the readers noting that case b) is similar
to that in Corollary~\ref{CON2}.

In case b) we have
$$\begin{vmatrix}
p_1& p_2& p_{i+1}\\
p_1'& p_2'& p_{i+1}'\\
x_1& x_2& x_{i+1}\\
\end{vmatrix}=0$$
and if $p_1:p_2:\dots\neq p_1':p_2':\dots,$ then
$$x_{i+1}=\frac{x_1(p_2p_{i+1}'-p_{i+1}p_2')+x_2(p_{i+1}p_1'-p_1p_{i+1})}
{p_2p_1'-p_1p_2'},\quad i\in\Bbb N.$$

Set $\tilde{x}_i=x_i/p_i,$ $i\in\Bbb N$ and
$$s_i=\frac{p_1+p_i}{p_1+p_2}+\frac{2}{p_1+p_2}\sum_{j=2}^{i-1}p_j,\quad i\ge 3.$$

\begin{cor}\label{CON4} If $p=p'\in \ell^1,$ then a sequence of
positive real numbers  $x=(x_i)_{i\in\Bbb N}$ is area attainable
if and only if $x_1>0,$
\begin{equation}\label{10}
\tilde{x}_i= \tilde{x}_1+s_i(\tilde{x}_2-\tilde{x}_1),\quad i\ge 3
\end{equation}
and
$$1-\frac{p_1+p_2}{p_1+2\sum_{j=2}^\infty p_j}<\frac{ \tilde{x}_2}{\tilde{x}_1}<2+\frac{p_2}{p_1}.$$

\end{cor}

\noindent {\bf Proof}. It is enough to check \eqref{10}. For that
purpose we rewrite \eqref{6} in the form
$$\frac{\tilde{x}_{j+1}-\tilde{x}_j}{p_{j+1}+p_j}=\frac{\tilde{x}_j-\tilde{x}_{j-1}}{p_j+p_{j-1}}.$$
Then
$$\tilde{x}_{j+1}-\tilde{x}_j=(p_{j+1}+p_j)\frac{\tilde{x}_2-\tilde{x}_1}{p_2+p_1}$$
and it remains to sum up these identities for
$j=2,\dots,i-1.$\qed

Note that for $p=p'\in\Bbb R^n,$ the above corollary is a more
compact expression of Proposition 4 in \cite{NN}.

{}

\begin{thebibliography}{}

\bibitem{NN} N. Nikolov, {\it On a problem from SMC} (in Bulgarian),
Matematika LXII (2023), no. 2, 30-33.

\end{thebibliography}
\end{document}